\documentclass{article}
\usepackage{amsmath}
\usepackage{amsthm}

\newcommand{\mrm}{\mathrm}
\newcommand{\mf}{\mathbf}

\newcommand{\ad}{{\mf{A}}}
\newcommand{\cc}{{\mf{C}}}
\newcommand{\htt}{{\mf{T}}}
\newcommand{\qq}{{\mf{Q}}}
\newcommand{\zz}{{\mf{Z}}}

\newcommand{\SL}{\mrm{SL}}
\newcommand{\GL}{\mrm{GL}}

\newcommand{\dnd}{{\not|}}

\newcommand{\expn}[1]{{e^{2\pi\sqrt{-1}#1}}}
\newcommand{\fdef}[3]{{{#1}\!:{#2}\rightarrow{#3}}}
\newcommand{\zmod}[1]{\zz/{#1}\zz}

\newlength{\setoftmpheight}
\newcommand{\setof}[2]{{\settoheight{\setoftmpheight}{${#2}$}\left\{
  \left.{#1}{\vrule width 0pt height \setoftmpheight}\ \right|{#2}\right\}}}

\newcommand{\mattwo}[4]{{\begin{pmatrix}{#1}&{#2}\\{#3}&{#4}\end{pmatrix}}}
\newcommand{\smallmattwo}[4]{{\left(\begin{smallmatrix}{#1}&{#2}\\
  {#3}&{#4}\end{smallmatrix}\right)}}

\newcommand{\MkGo}[1]{M_k(\Gamma_0({#1}))}
\newcommand{\MkGoN}{\MkGo{N}}
\newcommand{\MkGu}[1]{M_k(\Gamma^0({#1}))}
\newcommand{\MkGuN}{\MkGu{N}}

\newcommand{\SkGo}[1]{S_k(\Gamma_0({#1}))}
\newcommand{\SkGoN}{\SkGo{N}}
\newcommand{\SkGob}[1]{\overline{S}_k(\Gamma_0({#1}))}
\newcommand{\SkGoNb}{\SkGob{N}}

\newcommand{\ado}{{\ad^\infty}}
\newcommand{\GLtA}{{\GL_2(\ado)}}

\newcommand{\SLtN}{\SL_2(\zmod N)}
\newcommand{\SLtpi}{\SL_2(\zmod{p_i^{n_i}})}

\newcommand{\httN}{\htt^N}

\newtheorem{theorem}{Theorem}
\newtheorem{prop}[theorem]{Proposition}
\newtheorem{lemma}[theorem]{Lemma}

\title{On a Result of Atkin and Lehner}
\author{David Carlton}
\date{March 22, 1999}

\begin{document}

\maketitle

\section{Introduction}
\label{sec:intro}

We wish to give a new proof of one of the main results of
Atkin-Lehner~\cite{atkin-lehner}.  That paper depends, among other
things, on a slightly strengthened version of Theorem~\ref{thm:main}
below, which characterizes forms in $\SkGoN$ whose Fourier
coefficients satisfy a certain vanishing condition.  Our proof
involves rephrasing this vanishing condition in terms of
representation theory; this, together with an elementary linear
algebra argument, allows us to rewrite our problem as a collection of
local problems.  Furthermore, the classical phrasing of
Theorem~\ref{thm:main} makes the resulting local problems trivial;
this is in contrast to the method of Casselman~\cite{casselman}, whose
local problem relies upon knowledge of the structure of irreducible
representations of $\GL_2(\qq_p)$.  Our proof is therefore much more
accessible to mathematicians who aren't specialists in the
representation theory of $p$-adic groups; the method is also
applicable to other Atkin-Lehner-style problems, such as the level
structures that were considered in Carlton~\cite{carlton-thesis}.

Our proof of Theorem~\ref{thm:main} occupies Section~\ref{sec:main}.
In Section~\ref{sec:rest}, we explain the links between this Theorem
and the rest of Atkin-Lehner theory; in particular, we show that
Theorem~\ref{thm:main}, together with either the Global Result of
Casselman~\cite{casselman} or Theorem~4 of
Atkin-Lehner~\cite{atkin-lehner}, can be used to derive all of the
important results of Atkin-Lehner theory.

\section{The Main Theorem}
\label{sec:main}

Recall that, if $N | M$ and $d | (M/N)$, there is a map
$\fdef{i_d}{\MkGo{N}}{\MkGo{M}}$
defined by
\[
c_m(i_d(f)) = 
\begin{cases}
  0 & \text{if $d \dnd m$} \\
  c_{m/d}(f) & \text{if $d | m$}.
\end{cases}
\]
This map sends cusp forms to cusp forms and
eigenforms to eigenforms (with the same eigenvalues); up to
multiplication by a constant, it is given by $f \mapsto
f|_{\smallmattwo d001}$.

\begin{theorem}
  \label{thm:main}
  Let $f \in \MkGoN$ be such that $c_m(f) = 0$ unless $(m,N) > 1$.
  Then $f = \sum_{p | N} i_p(f_p)$, where $p$ varies over the primes
  dividing $N$ and where $f_p \in \MkGo{N/p}$.  Furthermore, if $f$ is
  a cusp form (resp.\ eigenform) then the $f_p$'s can be chosen to be
  cusp forms (resp.\ eigenforms with the same eigenvalues as $f$).
\end{theorem}

Our proof rests on two elementary linear algebra lemmas:

\begin{lemma}
  \label{lemma:la-ker}
  Let $V_1, \dots, V_n$ be vector spaces and, for each $i$, let $f_i$ 
  be an endomorphism of $V_i$.  Then
  \[
  \ker (f_1 \otimes \dots \otimes f_n) = \sum_{i=1}^n V_1 \otimes
  \dots \otimes (\ker f_i) \otimes \dots \otimes V_n.
  \]
\end{lemma}

\begin{proof}
  We can easily reduce to the case $n = 2$.  If we write $V_i = (\ker
  f_i) \oplus V'_i$ then $f_i|_{V'_i}$ is an isomorphism onto its
  image, and
  \[
  V_1 \otimes V_2 = ((\ker f_1) \otimes (\ker f_2)) \oplus ((\ker f_1)
  \otimes V'_2) \oplus (V'_1 \otimes (\ker f_2)) \oplus (V'_1 \otimes
  V'_2).
  \]
  We see that $f_1\otimes f_2$ kills the first three factors, and is
  an isomorphism from the fourth factor onto its image; $\ker(f_1
  \otimes f_2)$ is therefore the sum of the first three factors, which
  is what we wanted to show.
\end{proof}

\begin{lemma}
  \label{lemma:la-main}
  Let $V_1, \dots, V_n$ be vector spaces and, for each $i$, let
  $V'_i$ and $V''_i$ be subspaces of $V_i$.  Then
  \begin{multline*}
  \left(\sum_{i=1}^n V_1 \otimes \dots \otimes V'_i \otimes \dots \otimes
    V_n\right) \cap (V''_1 \otimes \dots \otimes V''_n) \\
  = \sum_{i=1}^n V''_1 \otimes \dots \otimes (V'_i \cap V''_i)
  \otimes \dots \otimes V''_n.
  \end{multline*}
\end{lemma}

\begin{proof}
  Again, we can assume that $n=2$.  Write $V_i = V_{i1} \oplus V_{i2}
  \oplus V_{i3} \oplus V_{i4}$ where $V_{i1} = V'_i \cap V''_i$, $V'_i
  = V_{i1} \oplus V_{i2}$, and $V''_i = V_{i1} \oplus V_{i3}$.  Then
  $V'_1 \otimes V_2 + V_1 \otimes V'_2$ is the direct sum of those
  $V_{1j}\otimes V_{2k}$'s where at least one of $j$ or $k$ is in the
  set $\{1,2\}$.  Also, $V''_1 \otimes V''_2$ is the direct sum of the
  $V_{1j}\otimes V_{2k}$'s where $j$ and $k$ are both in the set
  $\{1,3\}$.  Thus, their intersection is $(V_{11} \otimes V_{21})
  \oplus (V_{11} \otimes V_{23}) \oplus (V_{13} \otimes V_{21})$, as
  claimed.
\end{proof}

\begin{proof}[Proof of Theorem~\ref{thm:main}.]
  If $f \in \MkGoN$ then $f|_{\smallmattwo {N^{-1}}001} \in \MkGuN$,
  where we define the group $\Gamma^0(N)$ by
  \[
  \Gamma^0(N) = \setof{\mattwo abcd \in \SL_2(\zz)}{b\equiv 0 \pmod
    N}.
  \]
  Furthermore, up to multiplication by a constant, $f|_{\smallmattwo
    {N^{-1}}001}$ has the same Fourier coefficients as $f$, except
  that we have to take the $q$-expansion with respect to $\expn{z/N}$
  instead of $\expn{z}$.  Our Theorem, then, is equivalent to the
  statement that, if $f \in \MkGuN$ satisfies the condition
  \begin{equation}
    \label{eq:vanishing-coeff}
    \text{$c_m(f) = 0$ unless $(m,N)>1$}
  \end{equation}
  then $f = \sum_{p | N} f_p$ where $f_p \in \MkGu{N/p}$.

  Let $M = M_k(\Gamma(N))$; it comes with an action of $\SLtN$.
  If $f \in M$ and $d | N$, define $\pi_d(f)$ to be $\sum_{d|m}
  c_m(f)q^m$.  Then $\pi_d(f) \in M$: in fact,
  \[
  \pi_d(f) = \frac 1d\sum_{b=0}^{d-1} f|_\smallmattwo 1{bN/d}01.
  \]
  The principle of inclusion and exclusion implies that $f$ satisfies
  \eqref{eq:vanishing-coeff} if and only if
  \[
  f = \sum_{p | N} \pi_p(f) - \sum_{\substack{p_1,p_2|N\\p_1<p_2}}
  \pi_{p_1p_2}(f) + \dotsb.
  \]
  Thus, if $V$ is an irreducible $\SLtN$-representation contained in
  $M$, it suffices to prove our Theorem for a form in $V$, since the
  conditions of our Theorem can be expressed in terms of the action of
  $\SLtN$.
  
  Let $N = \prod_{i=1}^n p_i^{n_i}$ be the prime factorization of $N$.
  Then $\SLtN = \prod_i \SLtpi$, so $V = \bigotimes_i V_i$ where $V_i$
  is a representation of $\SLtpi$.  Also, $\pi_{p_i}$ acts as the
  identity on the $V_j$ for $j \ne i$.  So if we define
  \[
  \pi(f) = f - \sum_{p | N} \pi_p(f) +
  \sum_{\substack{p_1,p_2|N\\p_1<p_2}} \pi_{p_1p_2}(f) - \dotsb
  \]
  then $\pi = (1 - \pi_{p_1}) \otimes \dots \otimes (1 - \pi_{p_n})$
  and $\ker(\pi)$ is the space of forms satisfying
  \eqref{eq:vanishing-coeff}.  Thus, Lemma~\ref{lemma:la-ker} implies
  that
  \[
  \ker(\pi) = \sum_{i=1}^n V_1 \otimes \dots \otimes (\ker (1-\pi_{p_i}))
  \otimes \dots \otimes V_n.
  \]
  
  Turning now to the question of a form's being in $\MkGuN$, that is
  the case if and only if the form is both in $M_k(\Gamma(N))$ and is
  invariant under the image $B(N)$ of $\Gamma^0(N)$ in $\SLtN$.  Also,
  $B(N) = \prod_i B(p_i)$.  Thus, setting $V'_i$ to be
  $\ker(1-\pi_{p_i})$ and $V''_i$ to be the space of
  $B(p_i)$-invariant elements of $V_i$, Lemma~\ref{lemma:la-main}
  implies that an element of $V$ is both in $\ker \pi$ and invariant
  under $B(N)$ if and only if it is in
  \[
  \sum_{i=1}^n V''_1 \otimes \dots \otimes (V'_i \cap V''_i)
  \otimes \dots \otimes V''_n.
  \]
  But if $v_i \in V_i$ is in $V'_i \cap V''_i$ then it is invariant
  both under $B(p_i)$ and under projection to the subspace of
  invariants under the cyclic subgroup generated by $\smallmattwo
  1{p_i^{n_i-1}}01$; this last condition is equivalent to its being
  invariant under $\smallmattwo 1{p_i^{n_i-1}}01$.  Thus, our vector
  $v_i$ is invariant under
  \[
  \setof{\mattwo abcd \in \SLtpi}{b \equiv 0 \pmod{p_i^{n_i-1}}},
  \]
  and $V''_1 \otimes \dots \otimes (V'_i \cap V''_i) \otimes \dots
  \otimes V''_n$ is the set of invariants under
  \[
  \setof{\mattwo abcd \in \SLtN}{b \equiv 0 \pmod{N/p_i}},
  \]
  i.e. the elements of $V \cap \MkGu{N/p_i}$, completing our proof.
  
  The cusp form case is similar, replacing $M$ by the space of cusp
  forms.  The eigenform case then follows from the facts that the
  Hecke operators are simultaneously diagonalizable and that their
  action is preserved by the operators $i_p$.
\end{proof}

\section{Newforms, Oldforms, and All That}
\label{sec:rest}

In this Section, we explain the relation between
Theorem~\ref{thm:main} and the rest of Atkin-Lehner theory.  We shall
see that the whole theory follows from Theorem~\ref{thm:main} together
with facts about $L$-series associated to modular forms, as expressed
by Theorem~4 of Atkin-Lehner~\cite{atkin-lehner} or the Global Result
of Casselman~\cite{casselman}.  We claim no originality in the methods
used in this Section.

Define $K_0(N)$ to be the subspace of $f \in \SkGoN$ such that $c_m(f)
= 0$ unless $(m,N) > 1$: thus, $K_0(N)$ is the subspace characterized
in Theorem~\ref{thm:main}.  Define $\SkGoNb$ to be $\SkGoN/K_0(N)$;
for $f \in \SkGoNb$, $c_m(f)$ is well-defined exactly when $(m,N) =
1$.  Also, let $\httN$ be the free polynomial algebra over $\cc$
generated by commuting operators $T_m$ for $(m,N) = 1$.  Then $\httN$
acts on $\SkGoN$ (where $T_m$ acts as the $m$'th Hecke operator), and
its action is diagonalizable; it is easy to see that its action
descends to $\SkGoNb$.  (For example, $T_m$ commutes with the action
of the operators $\pi_d$ defined in the proof of
Theorem~\ref{thm:main}.)

\begin{prop}
  \label{prop:mult-one}
  The $\httN$-eigenspaces in $\SkGoNb$ are one-dimensional;
  furthermore, an eigenform $f \in \SkGoNb$ is zero if and only if
  $c_1(f) = 0$.
\end{prop}

\begin{proof}
  If $f \in \SkGoNb$ is an eigenform for $T_m$ with eigenvalue
  $\lambda_m(f)$ then $c_m(f) = \lambda_m(f)c_1(f)$.  Thus, if $f$ is
  a $\httN$-eigenform then it is determined by its eigenvalues and by
  $c_1(f)$.
\end{proof}

This Proposition, together with Theorem~\ref{thm:main}, sometimes
allows one to reduce questions about the spaces $\SkGoN$ to spaces
whose eigenspaces are one-dimensional.

\begin{prop}
  \label{prop:almost-all}
  If $f$ and $g$ are eigenforms in $\SkGoNb$ such that, for some $D$,
  they have the same eigenvalues $\lambda_m$ for all $m$ with $(m,ND)
  = 1$, then they have the same eigenvalues for all $m$ with
  $(m,N)=1$.
\end{prop}

\begin{proof}
  This is part of the Global Result of Casselman~\cite{casselman}, or
  of Theorem~4 of Atkin-Lehner~\cite{atkin-lehner}.
\end{proof}

We should also point out that our Theorem~\ref{thm:main} isn't quite
the same as Theorem~1 of Atkin-Lehner~\cite{atkin-lehner}.  Their
Theorem~1 assumes that $c_m(f) = 0$ unless $(m,ND) = 1$, and thus
breaks down into two parts: showing that you can assume that $D = 1$,
and our Theorem~1.  It is easy to show that the first part is
equivalent to Proposition~\ref{prop:almost-all}, at least in the
eigenform case; the cusp form case takes a bit more work.

We now present what is traditionally thought of as the core of
Atkin-Lehner theory.

\begin{theorem}
  \label{thm:growth}
  If $\{\lambda_m\}$ is a set of eigenvalues (for all $m$ relatively
  prime to a finite set of primes) that occurs in some space $\SkGoN$
  then there is a unique minimal such $N$ (with respect to division)
  for which those eigenvalues occur, and the corresponding eigenspace
  is one-dimensional.  If $f$ is a basis element for that eigenspace
  and if $M$ is a multiple of $N$ then the corresponding eigenspace in
  $\SkGo{M}$ has a basis given by the forms $i_d(f)$ where $d$ varies
  over the (positive) divisors of $M/N$.
\end{theorem}

\begin{proof}
  For any positive integer $M$, write $V_0(M)$ for the set of
  eigenforms in $\SkGo{M}$ with eigenvalues $\{\lambda_m\}$.  By
  Proposition~\ref{prop:almost-all}, we don't have to worry exactly
  about which primes are avoided in our set of eigenvalues, so this
  notation makes sense.  Furthermore, let $N$ be a minimal level such
  that $V_0(N)$ is nonzero.  By Proposition~\ref{prop:mult-one}, the
  image of $V_0(N)$ in $\SkGoNb$ is one-dimensional.
  Theorem~\ref{thm:main} shows that any element of the kernel of the
  map from $V_0(N)$ to $\SkGoNb$ is of the form $\sum_{p|N}i_p(f_p)$,
  where $f_p \in V_0(N/p)$.  But the minimality of $N$ shows that
  there aren't any such forms; the kernel is therefore zero, so
  $V_0(N)$ is one-dimensional.
  
  To see that $N$ is unique, let $S_k$ be the space of adelic cusp
  forms of weight $k$ but of arbitrary level structure; it comes with
  an action of $\GLtA$, and elements of $\SkGo{M}$ correspond to
  elements of $S_k$ invariant under the action of a certain subgroup
  $U_0(M) = \prod_p U_0(p^{m_p})$, where $p$ varies over the set of
  all primes and $p^{m_p}$ is the highest power of $p$ that divides
  $M$.  Casselman's Global Result says that the set $V$ of forms in
  $S_k$ with eigenvalues $\{\lambda_m\}$ gives an irreducible
  representation of $\GLtA$; thus, it can be written as a restricted
  tensor product $V = \bigotimes_p V_p$, and
  \[
  V_0(M) = \bigotimes_p V_p^{U_0(p^{m_p})}.
  \]
  Since $U_0(p^m)$ contains $U_0(p^{m+1})$, for each $p$ it is the
  case that, if for some power $p^m$, $V_p^{U_0(p^m)}$ is nonzero,
  then there is a minimal such power.  Thus, taking $N$ to be the
  product of those minimal powers of $p$, we see that, if for some
  $M$, $V_0(M)$ is nonzero, then it is nonzero for a unique minimal
  $M$, namely our $N$.  (Alternatively, the uniqueness of the minimal
  level is part of Theorem~4 of Atkin-Lehner~\cite{atkin-lehner}.)
  
  Finally, to see that the eigenspace grows as indicated, let $f$ be a
  nonzero element of $V_0(N)$ for $N$ minimal.  By
  Proposition~\ref{prop:mult-one}, we can assume that $c_1(f)=1$,
  since our argument above showed that the image of $f$ in $\SkGoNb$
  is nonzero.  Fix some multiple $M$ of $N$, and assume that we have
  shown that, for all proper divisors $M'$ of $M$ with $N | M'$,
  \begin{equation}
    \label{eq:direct-sum}
    V_0(M') = \bigoplus_{d|(M'/N)}i_d(f) \cdot \cc.
  \end{equation}
  We then want to show that the same statement holds with $M$ in place
  of $M'$.  Thus, let $g$ be an element of $V_0(M)$.  By
  Proposition~\ref{prop:mult-one}, the image of $g - c_1(g)i_1(f)$ in
  $\SkGob{M}$ is zero, so by Theorem~\ref{thm:main},
  \[
  g = c_1(g)i_1(f) + \sum_{p|M} i_p(g_p)
  \]
  for some forms $g_p \in V_0(M/p)$.  Also, $g_p = 0$ unless $p |
  (M/N)$, since otherwise $N$ wouldn't divide $M/p$, contradicting the
  unique minimality of $N$.  But then \eqref{eq:direct-sum} implies
  that each $g_p$, and hence $g$, can be written as a linear
  combination of the forms $i_d(f)$ for $d|(M/N)$; it is easy to see
  that such an expression for $g$ is unique.
\end{proof}

\end{document}